\documentclass[12pts,letter]{article}
\usepackage{hyperref}

\usepackage{multirow}
\usepackage{multicol} 

\usepackage[svgnames]{xcolor} 

\usepackage{times} 

\usepackage{graphicx} 
\usepackage{booktabs} 
\usepackage[font=small,labelfont=bf]{caption} 
\usepackage{amsfonts, amsmath, amsthm, amssymb} 
\usepackage[top=3cm, bottom=3cm, right=2cm, left=2cm]{geometry}

\title{{\color{DarkBlue} sPad of Super Accuracy and Geometry in Old Babylon.} 
\\sPad de S\'uper Exactitud y Geometr\'ia en la Antigua Babilonia.}
\author{F. Qui\~nonez\footnote{\texttt{faquinonez@saber.uis.edu.co}}, L. A. N\'u\~nez, F. D. Lora-Clavijo, 
\\R. Ort\'iz Aponte, O. Otero Olarte, D. Acosta Ort\'iz
\\\\Escuela de F\'isica - Maestr\'ia Matem\'atica Aplicada, 
\\Universidad Industrial de Santander, Bucaramanga, Colombia.}
\begin{document}
\maketitle
\begin{abstract}
{\color{DarkBlue} Recently it has been discovered that on a stone tablet over 3800 years old, the Plimpton-322 table,  
are carved the geometric relations that exist between the sides of 15
right triangles chosen in a very special way. Due to its property as a super accuracy calculation tool,
in this work we have called it sPad by stone pad, and we have calculated its machine accuracy $\epsilon_{m}$. 
Additionally, we present the physical and astrophysical constants most used in science and engineering in sexagesimal base.
}
\begin{center} ----- \end{center} 
Reci\'entemente se ha descubierto que en una tableta de piedra de m\'as de 3800
a\~nos de antig\"uedad, la Tabla Plimpton 322, est\'an ta\-lla\-das las relaciones 
geom\'etricas que existen entre los 
lados de 15 tri\'angulos rect\'angulos escogidos de manera muy especial. 
Debido a su propiedad como herramienta de c\'alculo s\'uper preciso,
en este trabajo la hemos llamado 
sPad por \emph{stone pad}, y hemos calculado su exactitud de m\'aquina $\epsilon_{m}$.
Adicionalmente presentamos las constantes f\'isicas y astrof\'isicas m\'as usadas
en ciencia e ingenier\'ia en base sexagesimal.
\end{abstract}
\section{Introducci\'on}
Los Antiguos Babilonios nos han dejado un registro hist\'orico de sus eventos mediante el uso de la escritura cuneiforme. 
Nues\-tro conocimiento de su matem\'atica proviene de unas 400 tabli\-llas grabadas en arcilla, 
las cuales abordan temas como fracciones, ecuaciones cuadr\'aticas, 
c\'ubicas y trios de enteros en aplicaci\'on del Teorema de Pit\'agoras (1000 a\~nos antes que Pit\'agoras).
sPad es conocida alrededor del mundo gracias al inter\'es que mostr\'o en ella Otto Neugebauer, 
Presidente del Departamento de Historia de Matem\'aticas de la Universidad de Brown. \'El y su asistente Sachs, 
interpretaron que dicha tabla conten\'ia algunas soluciones enteras del Teorema de Pit\'agoras \cite{PLIMPTON322}.
Her\'on de Alejandr\'ia, reconocido ingeniero y matem\'atico, considerado como uno de los cient\'ificos e inventores 
ma\'s grandes de la antiguedad, desarroll\'o t\'ecnicas de c\'alculo num\'erico tomadas de los Babilonios, 
tales como el c\'alculo de raices cuadradas mediante iteraciones. 
La f\'ormula de Her\'on establece la relaci\'on entre el \'area del tri\'angulo y sus tres lados.
$
A = \sqrt{s(s-a)(s-b)(s-c)}
$
Donde $a$, $b$ y $c$, son los lados de un tri\'angulo y $s$ su semiper\'imetro. 
El m\'etodo de Her\'on para el c\'alculo de raices, consiste en determinar los t\'erminos de la sucesi\'on definida por la recurrencia
$a_{n+1} = (a_{n} + a/a_{n} )/2$, 
tal que 
$\sqrt{a} = \lim_{n\rightarrow\infty} a_{n+1}$,
en donde $a$ es el n\'umero de la ra\'iz cuadrada que se desea calcular y $a_{0}$ una 
aproximaci\'on de la ra\'iz que buscamos. Es un m\'etodo que converge muy r\'apidamente \cite{ArticleP322}.
\begin{figure}
\centering
\includegraphics[width=13cm]{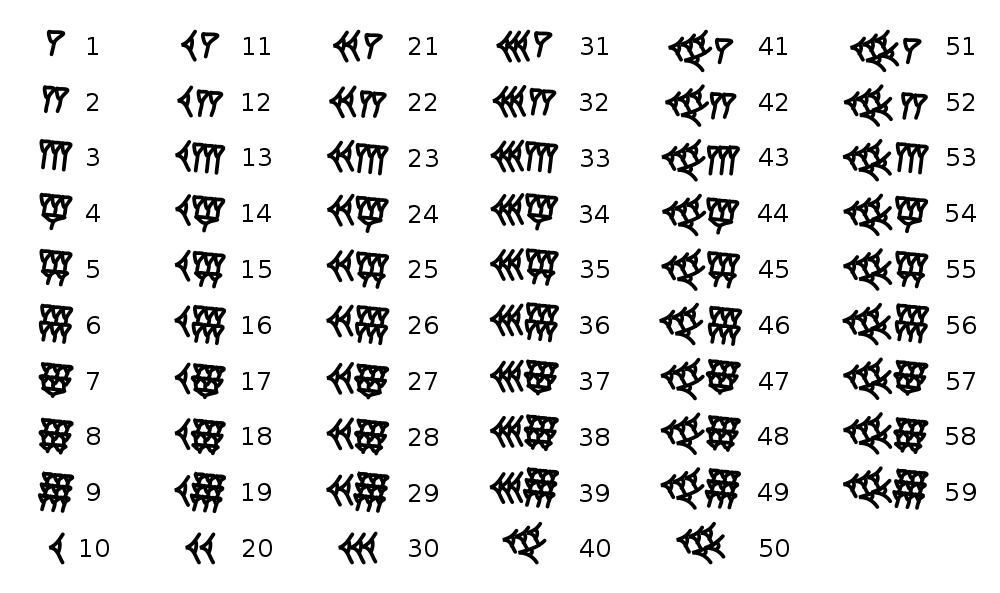}
\caption{Tomada de \url{https://en.wikipedia.org/wiki/File:Babylonian_numerals.svg}}
\label{f:numerals}
\end{figure}
\begin{figure}
\centering
\includegraphics[width=13cm]{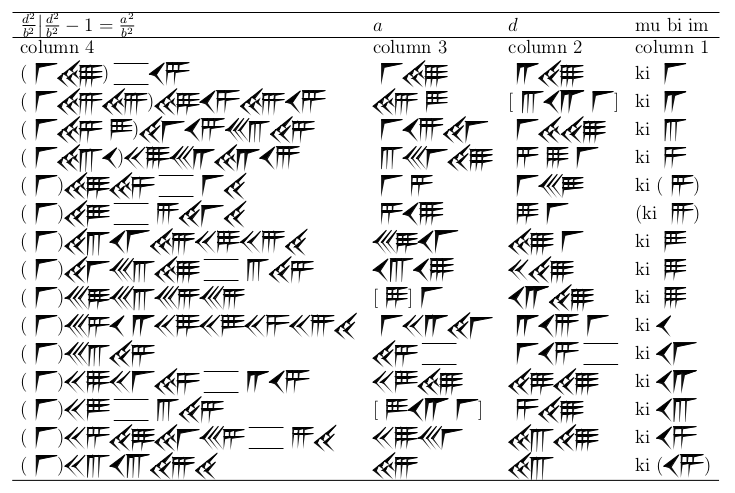}
\caption{Tabla Plimpton 322 llamada en este trabajo sPad \cite{PLIMPTON322}.}
\label{f:p322}
\end{figure}
\subsection{?`Por Qu\'e la Base Sexagesimal es mejor que la Decimal?}
Considerando descomposici\'on en n\'umeros primos
$10 = 2^{1} \times 5^{1},$
la base $10$ tiene $1+1=2$ divisiones enteras exactas: $10/2$ y $10/5$.
$60 = 2^{2} \times 3^{1} \times 5^{1}$,
la base $60$ tiene $2+1+1+2\binom{3}{2}=10$ divisiones enteras exactas:
$60/2$, $60/3$, $60/4$, $60/5$, 
$60/6$, $60/10$, $60/12$, $60/15$,
$60/20$, $60/30$.

\section{M\'etodos}
En la Figura~\ref{f:numerals} podemos ver los numerales sexagesimales como eran representados por los 
Antiguos Babilonios.
A continuaci\'on mostramos el mapa de caracteres usados en este trabajo.
\subsection{Mapa Decimal a Sexagesimal.}
		$0_{10}=\mathsf{0}$,
		$1_{10}=\mathsf{1}$,
		$2_{10}=\mathsf{2}$,
		$3_{10}=\mathsf{3}$,
		$4_{10}=\mathsf{4}$,
		$5_{10}=\mathsf{5}$,
		$6_{10}=\mathsf{6}$,
		$7_{10}=\mathsf{7}$,
		$8_{10}=\mathsf{8}$,
		$9_{10}=\mathsf{9}$,
		$10_{10}=\mathsf{A}$,
		$11_{10}=\mathsf{B}$,
		$12_{10}=\mathsf{C}$,\\
		$13_{10}=\mathsf{D}$,
		$14_{10}=\mathsf{E}$,
		$15_{10}=\mathsf{F}$,
		$16_{10}=\mathsf{G}$,
		$17_{10}=\mathsf{H}$,
		$18_{10}=\mathsf{I}$,
		$19_{10}=\mathsf{J}$,
		$20_{10}=\mathsf{K}$,
		$21_{10}=\mathsf{L}$,
		$22_{10}=\mathsf{M}$,
		$23_{10}=\mathsf{N}$, \\
		$24_{10}=\mathsf{O}$,
		$25_{10}=\mathsf{P}$,
		$26_{10}=\mathsf{Q}$,
		$27_{10}=\mathsf{R}$,
		$28_{10}=\mathsf{S}$,
		$29_{10}=\mathsf{T}$,
		$30_{10}=\mathsf{U}$,
		$31_{10}=\mathsf{V}$,
		$32_{10}=\mathsf{W}$,
		$33_{10}=\mathsf{X}$,
		$34_{10}=\mathsf{Y}$, \\
		$35_{10}=\mathsf{Z}$, 
		$36_{10}=\mathsf{\alpha}$,
		$37_{10}=\mathsf{\beta}$,
		$38_{10}=\mathsf{\gamma}$,
		$39_{10}=\mathsf{\delta}$,
		$40_{10}=\mathsf{\epsilon}$,
		$41_{10}=\mathsf{\zeta}$,
		$42_{10}=\mathsf{\eta}$, 
		$43_{10}=\mathsf{\theta}$,
		$44_{10}=\mathsf{\iota}$,
		$45_{10}=\mathsf{\kappa}$,
		$46_{10}=\mathsf{\lambda}$,
		$47_{10}=\mathsf{\mu}$,
		$48_{10}=\mathsf{\nu}$,
		$49_{10}=\mathsf{\xi}$,
		$50_{10}=\mathsf{o}$,
		$51_{10}=\mathsf{\pi}$,
		$52_{10}=\mathsf{\rho}$,
		$53_{10}=\mathsf{\sigma}$,
		$54_{10}=\mathsf{\tau}$,
		$55_{10}=\mathsf{\upsilon}$,
		$56_{10}=\mathsf{\phi}$,
		$57_{10}=\mathsf{\chi}$,
		$58_{10}=\mathsf{\psi}$,
		$59_{10}=\mathsf{\omega}$ 
y repite c\'iclicamente de la manera habitual $60_{10}=\mathsf{10}\ldots$.
En la Figura~\ref{f:p322} podemos ver la versi\'on original mientras que en la Tabla~\ref{t:p322}
vemos una versi\'on con s\'imbolos modernos que podemos usar con \LaTeX.
\subsection{Potencias de 60.}
$60^{-1}_{10}=0.01\overline{6}_{10}=\mathsf{10^{-1}}$ \\
$60^{-2}_{10}=0.0002\overline{7}_{10}=\mathsf{10^{-2}}$ \\
$60^{-3}_{10}=0.000004\overline{629}_{10}=\mathsf{10^{-3}}$ \\
$60^{-4}_{10}=0.00000007\overline{716049382}_{10}=\mathsf{10^{-4}}$ \\
$60^{-5}_{10}=0.00000000128600823045267496_{10}=\mathsf{10^{-5}}$ \\
$60^{-6}_{10}=0.0000000000214334705075445825_{10}=\mathsf{10^{-6}}$ \\
$60^{-7}_{10}=0.000000000000357224508459076342_{10}=\mathsf{10^{-7}}$ \\
$60^{-8}_{10}=0.00000000000000595374180765127242_{10}=\mathsf{10^{-8}}$ \\
$60^{-9}_{10}=0.00000000000000009.9229030127521216_{10}=\mathsf{10^{-9}}$
\begin{table}
\centering
\begin{tabular}{lrrc}
\toprule
$\frac{d^{2}}{b^{2}}$ o $\frac{a^{2}}{b^{2}}$ & $a$ & $d$ & fila \\
\midrule
$\mathsf{ 1 \omega 0 F }$                    & $\mathsf{ 1 \omega }$ & $\mathsf{ 2 \xi }$ & $\mathsf{ 1 }$ \\
$\mathsf{ 1 \varphi \varphi \psi E o 6 F }$        & $\mathsf{ \phi 7 }$   & $\mathsf{ 1 K P }$ & $\mathsf{ 2 }$ \\
$\mathsf{ 1 \upsilon 7 \zeta F X \kappa }$   & $\mathsf{ 1 G \zeta }$   & $\mathsf{ 1 o \xi }$ & $\mathsf{ 3 }$ \\
$\mathsf{ 1 \sigma A T W \rho G }$   & $\mathsf{ 3 V \xi }$   & $\mathsf{ 5 9 1 }$ & $\mathsf{ 4 }$ \\
$\mathsf{ 1 \nu \tau 1 \epsilon }$   & $\mathsf{ 1 5 }$   & $\mathsf{ 1 \beta }$ & $\mathsf{ 5 }$ \\
$\mathsf{ 1 \mu 6 \zeta \epsilon }$   & $\mathsf{ 5 J }$   & $\mathsf{ 8 1 }$ & $\mathsf{ 6 }$ \\
$\mathsf{ 1 \theta B \phi S Q \epsilon }$   & $\mathsf{ \gamma B }$   & $\mathsf{ \omega 1 }$ & $\mathsf{ 7 }$ \\
$\mathsf{ 1 \zeta X \kappa E 3 \kappa }$   & $\mathsf{ D J }$   & $\mathsf{ K \xi }$ & $\mathsf{ 8 }$ \\
$\mathsf{ 1 \gamma X \alpha \alpha }$   & $\mathsf{ 8 1 }$   & $\mathsf{ C \xi }$ & $\mathsf{ 9 }$ \\
$\mathsf{ 1 Z A 2 S R O Q \epsilon }$   & $\mathsf{ 1 M \zeta }$   & $\mathsf{ 2 G 1 }$ & $\mathsf{ A }$ \\
$\mathsf{ 1 X \kappa }$   & $\mathsf{ \kappa }$   & $\mathsf{ 1 F }$ & $\mathsf{ B }$ \\
$\mathsf{ 1 T L \tau 2 F }$   & $\mathsf{ R \omega }$   & $\mathsf{ \nu \xi }$ & $\mathsf{ C }$ \\
$\mathsf{ 1 R 0 3 \kappa }$   & $\mathsf{ 2 \zeta }$   & $\mathsf{ 4 \xi }$ & $\mathsf{ D }$ \\
$\mathsf{ 1 P \nu \pi Z 6 \epsilon }$   & $\mathsf{ T V }$   & $\mathsf{ \sigma \xi }$ & $\mathsf{ E }$ \\
$\mathsf{ 1 N D \lambda \epsilon }$   & $\mathsf{ \phi }$   & $\mathsf{ 1 \lambda }$ & $\mathsf{ F }$ \\
\bottomrule
\end{tabular}
\caption{Representaci\'on moderna de la Tabla Plimpton 322. La primer columna presenta una ambig\"uedad en el primer sexag\'esito
puesto que no se tiene certeza en que valor de significancia inicial, los Babilonios representaban los n\'umeros flotantes en base
sexagesimal. Las columnas 2 y 3 en la tabla original vienen justificadas a la izquierda, $a$ y $d$ representan el lado m\'as corto y la diagonal de un 
tri\'angulo rect\'angulo, respectivamente, a su vez $b$ es el lado mediano.}
\label{t:p322}
\end{table}
%
%
\subsection{N\'umeros Flotantes.}
Sea $a\in \mathbb{R}$ en las m\'aquinas modernas de c\'omputo se representan los n\'umeros de punto flotante mediante
\begin{equation}
a = s \cdot M \cdot B^{e-E}
\end{equation}
en donde $B$ es la base que en su mayor\'ia es 2, $s$ es un bit que se deja para el signo, 
$e$ es el exponent, $E$ el bias del exponente y $M$ es la mantissa que viene normalizada en la 
mayor\'ia de los casos $M \in [\frac{1}{B},1)$.
Nosotros apreciamos que la primera columna de la sPad es una mantissa normalizada en base sexagesimal
que no usa exponente en la base. Creemos que la raz\'on de no usar exponente es sencillamente 
porque ellos no necesitaban m\'as exactitud en sus c\'alculos y tampoco necesitaban representar n\'umeros
muy peque\~nos o muy grandes, por lo tanto
\begin{equation}
a = M 
\end{equation}
de manera tal que en base sexagesimal 
\begin{equation}
M = 
\begin{cases}
m_{1} \mathsf{10^{-1}} +  \ldots + m_{8} \mathsf{10^{-8}} \\
m_{1} \mathsf{10^{0}} + \ldots + m_{9} \mathsf{10^{-8}} \\
\end{cases}
\end{equation}
preservando la ambig\"uedad si representa a $s^{2}/l^{2}$ o $d^{2}/l^{2}$.
Las columnas 2, 3 y 4 representan n\'umeros en base sexagesimal de la manera habitual:
\begin{equation}
a = a_{u-1} \mathsf{10^{u-1}} + \ldots + a_{1} \mathsf{10^{1}} + a_{0}
\end{equation}
en donde reconocemos el real $a$ se escribe mediante $u$ sexag\'esitos como $a_{u-1}\ldots a_{1}a_{0}$.
\begin{equation} 
\begin{aligned}
10^{-n} &= 6^{n} \; 60^{-n} \\
10^{n} &= 6^{-n} \; 60^{n}
\end{aligned}
\end{equation}
Hemos usado \cite{ROOT} en este trabajo para realizar los c\'alculos concernientes a grandes enteros de 8 bytes.
%
\newpage
\section{Resultados}
La exactitud de m\'aquina de la sPad es $60^{-8}$.
\begin{equation}
\epsilon_{m} = 60^{-8} = 5.95374180765127242 \times 10^{-15} 
\end{equation}
\begin{table}[h!]
\centering
\begin{tabular}{ccrl}
\toprule
        Constante & S\'imbolo & Valor Sexagesimal &  \textrm{Unidad} \\
		\midrule
        Const. de Avogrado & $N_A$ & $\mathsf{4\alpha\delta BI\rho WCK\varphi 4}$ & [mol$^{-1}$] \\
		Velocidad de la luz & $c$ & $\mathsf{N7\nu\epsilon\psi}$ &  $[ms^{-1}]$ \\
        Const. Planck & $h/2\pi$ & $\mathsf{1P\psi 1MYM\xi 3A\cdot 10^{-L} }$ & [eV$\cdot$s]\\
        Mag. carga electron& $e$ & $\mathsf{5\nu \kappa WI\psi P\omega\eta \nu \cdot 10^{-K} }$ & [C]\\
        Const. de conversi\'on & $(\hbar c)$ & $\mathsf{3HJ\beta7PE\rho \nu \cdot 10^{-7}}$ & [MeV*fm]\\
        Const. de conversi\'on & $(\hbar c)^2$ & $\mathsf{NL\kappa\varphi Y\zeta Q0\psi \cdot 10^{-9}}$ & [GeV$^2\cdot$mbarn]\\
        Masa electr\'on & $m_e$ & $\mathsf{U\delta Z\lambda KTET\kappa\cdot 10^{-9} }$ & [MeV/c$^2$]\\
        Masa prot\'on & $m_p$ & $\mathsf{ FyGJTXy\rho \nu\cdot 10^{-7} }$ & [MeV/c$^2$]\\
        Const. gravitacional & $G_N$ & $\mathsf{36\xi\sigma TZUUV\alpha\cdot 10^{-F} }$ & [m$^3\cdot$s$^{−2}$/kg] \\
        Ac. gra. estandar& $g_N$ & $\mathsf{ 9\nu N\varphi O0\cdot 10^{-5} }$ & [m/s$^2$] \\
        Const. de Boltzmann & $k$ & $\mathsf{1vB\alpha 7I\delta IRK\cdot 10^{-M} }$ & [J/K] \\
        Masa deuter\'on & $m_d$ & $\mathsf{ VF\alpha\lambda WQ\rho\nu \cdot 10^{-6} }$ & [MeV/c$^2$] \\
        Unidad masa at\'omica & $1g/N_A mol$ & $\mathsf{ FVTy\iota\alpha N2O\cdot 10^{-7} }$ & [MeV/c$^2$] \\
\bottomrule
\end{tabular}	
\label{tab:Plimpton322}
\caption{Constantes f\'isicas en base sexagesimal. Constantes tomadas de \cite{PDG2017}.}
\end{table}

\section{Conclusiones}
Los antiguos Babilonios:
\begin{itemize}
\item Usaban base num\'erica sexagesimal como sistema de numeraci\'on.
\item Ten\'ian la habilidad de calcular relaciones entre los lados de cualquier tri\'angulo rect\'angulo con una exactitud de al menos 15 cifras decimales.
\item Usaban el concepto de mantissa normalizada en el modelado
de la sPad. F\'acilmente pod\'ia extenderse la exactitud en sus sPads ya que s\'olo depende del tama\~no f\'isico de la tabla para seguir agregando
mas sexag\'esitos de significancia a la mantissa tallada, aunque no sabemos si pod\'ian usar dicha exactitud en la pr\'actica.
\item Fabricaron sPads que han demostrado ser unos dispositivos portables de larga durabilidad y no requiere del uso de baterias o combustible f\'osil para su funcionamiento.
Tambi\'en pod\'ian ser usadas por personas con discapacidad visual que manejaran sensibilidad en sus manos para leer los tallados.
\end{itemize}



\bibliographystyle{plain} 
\bibliography{sexagesimal} 



\end{document}